\documentclass[10pt]{article}

\usepackage{amsmath}
\usepackage{amssymb} 
\usepackage{latexsym} 
\usepackage{epic,eepic,epsfig,psfrag} 
\usepackage{amscd}  
\usepackage{index}  

\newcommand{\QED}{\hspace*{\fill}$\Box$\medskip} 

\def\one{\hbox{1\hskip-2.7pt l}}

\def\rd{{\rm d}}

\def\rG{{\rm G}}
\def\p{\phi}
\def\ph{\varphi}

\def\k{\kappa}
\def\ep{\varepsilon} 
\def\e{\eta}

\def\x{\xi} 
 
\def\n{\nu}

\def\pg{\mathhexbox278}
\def\S{\Sigma}

\def\X{\Xi} 
\def\Om{\Omega} 
\def\P{\Phi} 
\def\cg{\mathfrak{g}}
\def\cA{{\mathcal A}}

\def\cE{{\mathcal E}} 
 
\def\cG{{\mathcal G}}

\def\CS{{\mathcal C}{\mathcal S}}
\def\E{{\mathcal E}}
\def\R{{\mathbb R}} 
 
\def\T{{\mathbb T}} 
\def\N{{\mathbb N}} 
\def\C{{\mathbb C}} 
\def\CP{{\mathbb C\mathbb P}} 
 
\def\Z{{\mathbb Z}}
 
\def\half{{\textstyle{\frac 12}}}

\def\Pr{{\bf Proof:}\;}

\def\dph{{\rm d}\p}

\def\dr{{\rm d}r}
\def\pd{\partial}

\def\la{\langle\,}
\def\ra{\,\rangle}

\def\SU{{\rm SU}}
\def\su{{\mathfrak s\mathfrak u}}
\def\SO{{\rm SO}}
\def\so{{\mathfrak s\mathfrak o}}

\def\tint{\textstyle\int}

\newtheorem{dfn}{Definition}[section] 
\newtheorem{lem}[dfn]{Lemma} 
 
\newtheorem{thm}[dfn]{Theorem} 
\newtheorem{rmk}[dfn]{Remark}

\def\co{\colon\thinspace}

\begin{document}

\bibliographystyle{plain}

\author{Katrin Wehrheim}

\title{Energy identity for anti-self-dual instantons on $\C\times\S$}




\maketitle

\begin{abstract}

We prove an energy identity for anti-self-dual connections on the product 
${\C\times\S}$ of the complex plane and a Riemann surface.
The energy is a multiple of a basic constant that is determined 
from the values of a corresponding Chern-Simons functional 
on flat connections and its ambiguity under gauge transformations.

For $\SU(2)$-bundles this identity supports the conjecture that the finite
energy anti-self-dual instantons correspond to holomorphic bundles
over $\CP^1\times\S$.

Such anti-self-dual instantons on $\SU(n)$- and $\SO(3)$-bundles arise in 
particular as bubbles in adiabatic limits occurring in the context of
mirror symmetry and the Atiyah-Floer conjecture.
Our identity proves a quantization of the energy of these bubbles that
simplifies and strengthens the involved analysis considerably.

\end{abstract}

\section{Introduction}

Let $\S$ be a Riemann surface and consider the trivial $\SU(2)$-bundle 
over $\C\times\S$.
A connection $\X\in\cA(\C\times\S)$ on this bundle is a $1$-form
$\X\in\Om^1(\C\times\S;\su(2))$ with values in the Lie algebra $\su(2)$.
Gauge transformations $u\in\cG(\C\times\S)$ of the bundle are
represented by maps $u\in{\rm Map}(\C\times\S,\SU(2))$ and act
on $\cA(\C\times\S)$ by $u^*\X=u^{-1}\X u + u^{-1}\rd u$.
We equip $\C\times\S$ with a product metric of the Euclidean metric on $\C$ and
a fixed metric on $\S$.
Then a connection $\X\in\cA(\C\times\S)$ is called an {\bf ASD instanton} if 
its curvature is anti-self-dual,
$$
F_\X + * F_\X = 0 ,
$$
where $*$ is the Hodge operator w.r.t.\ the metric on $\C\times\S$.
The curvature $F_\X=\rd\X + \X\wedge\X$ transforms under gauge transformations
$u\in\cG(\C\times\S)$ as $F_{u^*\X}=u^{-1}F_\X u$, hence the anti-self-duality
equation is gauge invariant.
Next, we equip $\su(2)$ with the $\SU(2)$-invariant inner product 
$\la\x,\e\ra=-{\rm tr}(\x\e)$.
Then the energy of a connection $\X\in\cA(\C\times\S)$ is the gauge
invariant quantity
\begin{align*}
\E(\X) := \half \tint_{\C\times\S} |F_\X|^2 .
\end{align*}

The main purpose of this note is to establish the following energy identity.
Its surprisingly simple proof is given in section~\ref{sec:energy}. 
For the sake of simplicity we first focus our attention to $\SU(2)$-bundles.
Later, we will also indicate how to generalize this result to other 
structure groups and nontrivial bundles over $\S$.

\begin{thm} \label{thm energy}
Let $\X\in\cA(\C\times\S)$ be an ASD instanton.
If it has finite energy \hbox{$\E(\X)<\infty$}, then actually
$\E(\X)\in 4\pi^2 \N_0$.
\end{thm}

This energy quantization supports a conjectural correspondence between finite energy 
ASD instantons on $\C\times\S$ and holomorphic bundles over $\CP^1\times\S$.
For $\S=\T^2$ 
Biquard and Jardim \cite{BJ} showed that the gauge equivalence classes of 
ASD instantons with quadratic curvature decay are in one-to-one correspondence 
to a class of rank 2 stable holomorphic bundles over $\CP^1\times\T^2$.
Here the holomorphic structure induced by an instanton $\X$ extends over 
$\{\infty\}\times\T^2$ to define a bundle $E$, whose second Chern number
is given by the instanton energy, 
$c_2(E)=\frac 1{8\pi^2}\int \la F_\X\wedge F_\X \ra$, see \cite[\pg 2.3]{J}.
By our result this formula continues to give integer (Chern ?) numbers 
for finite energy instantons and any surface~$\S$.

\begin{rmk} \label{rmk energy}
A similar energy identity holds for ASD instantons on $\C\times P$ 
for any principal bundle $P\to\S$ with compact structure group $\rG$: 

Suppose that the Lie algebra $\cg$ is equipped with a $\rG$-invariant metric 
that satisfies {\rm (H)} below.
Then the statement of theorem~\ref{thm energy} holds with $4\pi^2$ replaced 
by the constant $\k_\cg N_\rG^{-1}$ given below.
\end{rmk}

On a nontrivial bundle $P$ the gauge transformations are represented by 
sections in the associated bundle $\rG_P = P\times_c \rG$ 
(using the conjugation action on $\rG$).
We can pick a $\rG$-invariant inner product on $\cg$ (and thus on
$\cg_P=P\times_{\rm Ad} \cg$).
Then the Maurer-Cartan $3$-form on each fibre of $\rG_P$ induces a closed $3$-form
$\e_\rG := \frac 1{12} \la g^{-1}\rd g \wedge [ g^{-1}\rd g \wedge g^{-1}\rd g ] \ra$
on $\rG_P$.
We need the following assumption.\\


\noindent
{\bf (H):}
There exists $\k_\cg>0$ such that
$[\k_\cg \e_\rG] \in {\rm H}^3(\rG_P,\R)$
is an integral class.\\


This holds for example with $\k_{\so(3)}=4\pi^2$ for any $\SO(3)$-bundle 
when we choose the inner product $-2{\rm tr}(\x\e)$ for $\x,\e\in\so(3)$.
It can also be achieved for any simply connected compact Lie group $\rG$,\footnote{
In that case the bundle is automatically trivial and
$\rG$ is isomorphic to a product $S_1\times\ldots\times S_k$ of simply connected, 
simple, and compact Lie groups $S_j$ with $\pi_3(S_j)\cong\Z$.
So we can pick a metric on each factor $S_j$ for which $[\e_{S_j}]\in {\rm H}^3(S_j,\R)$
is integral.}
e.g.\ for the trivial $\SU(n)$-bundles.
Finally, $N_\rG$ is the least common multiple of $\{1,2,\ldots,n_\rG\}$,
where $n_\rG$ denotes the maximal number of connected components that the
centralizer of a subgroup in $\rG$ can have.
This is finite since $\rG$ is compact. 
For $\SO(3)$ we have ${N_{\SO(3)}=1}$.

One source of interest in the ASD instantons on $\C\times\S$ is
the following adiabatic limit.
Let $\S\hookrightarrow X\to M$ be a fibre bundle with $\dim X=4$.
Consider ASD instantons $\X_\ep$ over $X$ with respect to metrics $g_M + \ep^2 g_\S$ 
for a sequence $\ep\to 0$.
If $|F_{\X_\ep}|_\text{fibre}|+\ep^2|F_{\X_\ep}|_\text{mix}|$ converges to a nonzero value,
then local rescaling on $M$ (but not in the fibre) yields an ASD instanton on $\C\times\S$
in the limit. This bubbling phenomenon is a central difficulty of the limiting process.
Adiabatic limits of this type have fascinating consequences from topology to
mathematical physics. They were first considered by Dostoglou-Salamon \cite{DS}, and recently
by Chen \cite{C} and Nishinou \cite{Ni}.
The energy quantization presented here simplifies and strengthens the bubbling analysis
and results in all these cases.
It can also be used for the Atiyah-Floer conjecture project \cite{Sa, W survey}. \\

{\it I would like to thank Benoit Charbonneau, 
Kenji Fukaya, Marcos Jardim, Tom Mrowka, and Dietmar Salamon 
for helpful discussions.}

\section{Proof of the energy identity}
\label{sec:energy}

In the following, $S_r\subset\C$ denotes the circle of radius $r$ centered at $0$.
We moreover denote by $D_r\subset\C$ the disk of radius $r$, and we introduce 
polar coordinates $(r,\p)\in(0,\infty)\times S^1$ on $\C^*=\C\setminus\{0\}$, 
with $S^1=\R/2\pi\Z$. 
Then on $\C^*\times\S$ we can write a connection $\X\in\cA(\C\times\S)$ 
in the splitting
$$
\X = A(r) + R(r)\dr + \P(r) \dph
$$
with $A(r)\co S^1\to\cA(\S)$ and $R(r),\P(r)\co S^1\to\Om^0(\S,\su(2))$ for
all $r\in(0,\infty)$.
The anti-self-duality equation becomes in this splitting
\[
\left\{\begin{aligned} 
 r^{-1} \bigl(  \pd_r\P - \pd_\p R + [\P,R] \bigr) + *F_A &= 0 ,\\
r^{-1} \bigl( \pd_\p A - \rd_A\P  \bigr) - *\bigr( \pd_r A - \rd_A R \bigr) &= 0 .
\end{aligned}\right.
\]
By $F_\X(r)$ we denote the curvature of $\X\in\cA(\C\times\S)$
over $S_r\times\S$ (but as a $2$-form on $\C\times\S$).
Then the curvature of an ASD instanton is
$$
\half |F_\X(r)|^2 = | F_{A(r)} |^2 + r^{-2} \bigl| \pd_\p A(r) - \rd_{A(r)}\P(r) \bigr|^2 .
$$
The energy of an ASD instanton on $D_r\times\S$ can be 
expressed in terms of the Chern-Simons functional of
$B(r) := A(r) + \P(r) \dph \in\cA(S^1\times\S)$,
$$
\half \int_{D_r\times\S} |F_\X|^2 
\;=\; - \half \int_{D_r\times\S} \la F_\X \wedge F_\X \ra 
\;=\; -\CS(B(r)) .
$$
The Chern-Simons functional on connections $B=A+\P\dph \in \cA(S^1\times P)$ 
is 
\begin{align}
\CS(B) 
&= \half \int_{S^1\times\S} \la B \wedge \bigl( F_B - \tfrac 16 [B\wedge B] \bigr) \ra  \nonumber\\
&= \int_{S^1}\int_\S \half \la \pd_\ph A \wedge A \ra + \la F_A \,,\, \P \ra . \label{CS}
\end{align}
For future reference we note the following identity which shows that the Chern-Simons functional
is continuous with respect to the $W^{1,\frac 32}$-norm. (Note that 
$W^{1,\frac 32}\hookrightarrow L^3$ on a $3$-manifold.)
For all $B,B^0\in\cA(S^1\times\S)$
\begin{align}
&\CS(B) - \CS(B^0)  \label{CS diff}\\
&=  \int \half \la ( F_B + F_{B^0}) \wedge (B- B^0) \ra      
        -\tfrac 1{12} \la [(B-B^0)\wedge(B-B^0)]\wedge (B-B^0) \ra.  \nonumber
\end{align}
The Chern-Simons functional is not gauge invariant, but its ambiguity on gauge orbits
is determined by the degree of the gauge transformations (as maps to $\SU(2)\cong S^3$): 
For all $B\in\cA(S^1\times\S)$ and $u\in\cG(S^1\times\S)$
\begin{equation} \label{CS deg}
\CS(B) - \CS(u^*B) = 4\pi^2 {\rm deg}(u) \in 4\pi^2\,\Z
\end{equation} 
For a general (possibly nontrivial) bundle $P\to\S$ one has to fix a flat reference connection. 
Then connections are given by $1$-forms with values in $\cg_P$ and 
the Chern-Simons functional depends on the choice of this reference 
connection only up to an additive constant.
(The proof of theorem~\ref{thm energy} will show that a flat connection exists.)
The right hand side of (\ref{CS diff}) is then given by $\int u^* \e_\rG$. 
So under the assumption {\rm (H)} we have $\CS(B) - \CS(u^*B) \in \k_\cg\,\Z$.

The second point that affects the constant in the energy identity is the possible values 
of the Chern-Simons functional on flat connections. The following result holds for
$\SU(2)$- and $\SO(3)$-bundles, and we will give the argument for a general bundle 
$P\to\S$, indicating how to proceed for other structure groups.

\begin{lem} \label{lem flat}
For every flat connection $B\in\cA_{\rm flat}(S^1\times\S)$ there is a gauge transformation 
$u\in\cG(S^1\times\S)$ such that $\CS(u^*B)=0$, and consequentially
$\CS(B)= 4\pi^2{\rm deg}(u)\in 4\pi^2\,\Z$.
\end{lem}

\Pr
We periodically extend the given flat connection to $B\in\cA_{\rm flat}(\R\times P)$. Then
one can find $u\co \R\to\cG(P)$ such that $u(0)\equiv\one$ and $u^*B\in\cA_{\rm flat}(\R\times P)$
has no $\dph$-component. Thus the curvature component $\pd_\p(u^*B)$ vanishes, 
and hence $u^*B\equiv A^0 \in \cA_{\rm flat}(P)$.
This is done by solving $\pd_\p u = -\P u$, so due to the periodicity of $\P$ we obtain 
the twisted periodicity $u(\p+2\pi)=u(\p)u(2\pi)$ for the gauge transformation.
Unless $u(2\pi)\equiv\one$ this does not define a gauge transformation on $S^1\times P$.
However, we know that $u(2\pi)$ lies in the isotropy subgroup $\cG_{A^0}$, since
$u(2\pi)^*A^0=u(2\pi)^*B(2\pi,\cdot)|_\S=u(0)^*B(0,\cdot)|_\S=A^0$.
If $\cG_{A^0}$ is connected, then we can multiply $u$ with a path within $\cG_{A^0}$ 
from $\one$ to $u(2\pi)^{-1}$ to obtain the required gauge transformation 
$w\in\cG(S^1\times P)$.
It satisfies $w^*B= A^0 + \P^0\dph$ with $\pd_\p A^0=0$ but possibly nonzero $\P^0$.
Now compare (\ref{CS}) to see that $\CS(w^*B)=0$,
and so $\CS(B)= 4\pi^2{\rm deg}(w)$ by (\ref{CS deg}).

For $\SO(3)$-bundles, any isotropy subgroup is connected since 
any centralizer (of the holonomy subgroup) in $\SO(3)$ is connected. 
Thus the proof is finished.
For a general Lie group whose centralizers have up to $n_\rG$ components, 
one finds that $u(2\pi)^n$ is homotopic to the identity for some integer $n\leq n_\rG$.
Then an "$n$-fold cover" $B^{(n)}$ of $B$ can be put into a gauge whose Chern-Simons
functional vanishes, and thus 
$\CS(B) = n^{-1} \CS(B^{(n)}) \in \k_\cg n^{-1}\,\Z \subset \k_\cg N_\rG^{-1}\,\Z$
if {\rm (H)} holds.

For $\SU(2)$ we would have $n_\rG=2$ due to the centralizer $\{\one,-\one\}$.
However, since the isotropy element $u(2\pi)=-\one$ is a constant, 
we do not need to go to a cover.
More generally suppose that $u(2\pi)=\exp(2\pi\x)$ for some constant $\x\in\cg$.
Let $v(\p):=\exp(-\p\x)$, then $w:=uv\in\cG(S^1\times\S)$ and
$w^*B= v^{-1}A^0 v - \x\dph$ (and both are of class $W^{1,\infty}$).
Then using $F_{A^0}=0$ and $\rd\x=0$ we obtain
\begin{align*}
\CS(w^*B) 
&= \int_{S^1}\int_\S \half \la v^{-1} [\x,A^0] v \wedge v^{-1}A^0 v \ra \\
&= \int_{S^1}\int_\S \la \x , A^0 \wedge A^0 \ra 
\;=\; -  \int_{S^1}\int_\S \la \x , \rd A^0 \ra \;=\; 0.
\end{align*}
\QED

In the subsequent proof of the energy identity we also work with a general
bundle $P\to\S$ and only for the final conclusion use the knowledge
from lemma~\ref{lem flat} on the 
possible values of the Chern-Simons functional on flat connections. \\ 

\noindent
{\bf Proof of theorem \ref{thm energy}: }
Note that for $B(r)\in\cA(S^1\times\S)$ given by $\X$ on $S_r\times\S$ we have
$|F_{B(r)}|^2 = |F_{A(r)}|^2 + |\pd_\ph A(r) - \rd_{A(r)}\P(r)|^2 
\leq \half r^2 |F_\X(r)|^2$
when $r\geq 1$, so
$$
\int_1^\infty r^{-1} \|F_{B(r)}\|_{L^2(S^1\times\S)}^2 \dr \;\leq\; \E(\X) < \infty .
$$
Thus we find a sequence $r_i\to\infty$ with $\|F_{B(r_i)}\|_{L^2(S^1\times\S)}\to 0$.
By Uhlenbeck's weak compactness \cite{U} we then find a further subsequence,
gauge transformations $u_i\in\cG(S^1\times P)$, and a flat limit connection 
$B^\infty\in\cA_{\rm flat}(S^1\times P)$ such that
$\| u_i^*B(r_i) - B^\infty\|_{W^{1,2}(S^1\times\S)}\to 0$.\footnote{
Originally, this convergence is only in the weak $W^{1,2}$-topology and in the
$L^4$-norm due to a compact Sobolev embedding.
By a local slice theorem (e.g.\ \cite[Theorem~8.1]{W}) one can 
achieve the additional relative Coulomb gauge condition
$\rd_{A_\infty}^*(v_\n^*A_\n - A_\infty)=0$.
Moreover, we have
$\rd_{A_\infty}(v_\n^*A_\n - A_\infty) 
=v_\n^{-1}F_{A_\n}v_\n - \half [(v_\n^*A_\n - A_\infty)\wedge (v_\n^*A_\n - A_\infty)] $.
So from the regularity of the Hodge decomposition of $1$-forms (e.g.\ \cite[Theorem~5.1]{W})
one obtains the convergence in the $W^{1,2}$-norm.
}
So we have $\CS(u_i^*B(r_i)) \to \CS(B^\infty)$ due to the convergence of $u_i^*B(r_i)$
and (\ref{CS diff}).
On the other hand the energy is finite, so
\begin{align*}
\E(\X) \;=\; - \lim_{i\to\infty} \int_{D_{r_i}\times\S} \la F_\X \wedge F_\X \ra  
\;=\; - \lim_{i\to\infty} \CS(B(r_i)).
\end{align*}
This shows that $\CS(B(r_i))$ also converges.
Now for an $\SU(2)$-bundle we have
$\CS(B^\infty)\in 4\pi^2\,\Z$ from  lemma~\ref{lem flat}.
Thus $\CS(B(r_i))= \CS(u_i^* B(r_i)) + 4\pi^2{\rm deg}(u_i)$ must converge
to some value in $4\pi^2\,\Z$. This proves the claim since that limit is also
the energy $\cE(\X)$.
%
%
For a general bundle under the assumption {\rm (H)} we
similarly obtain $\cE(\X) \in \k_\cg N_\rG^{-1}\,\Z$.
\QED

\pagebreak

\bibliographystyle{plain}

\end{document}